\documentstyle[12pt]{article}

\font\twelvemsb=msbm10 scaled 1200	
\font\eightmsb=msbm8			
\font\sixmsb=msbm6			
\newfam\msbfam
 \textfont\msbfam=\twelvemsb
 \scriptfont\msbfam=\eightmsb
 \scriptscriptfont\msbfam=\sixmsb
\def\Bbb#1{\fam\msbfam\relax#1}
\def\N{\Bbb{N}}
\def\Z{\Bbb{Z}}
\font\twelveeus=eusm10 scaled 1200	
\font\eighteus=eusm8			
\font\sixeus=eusm6			
\newfam\eusfam
 \textfont\eusfam=\twelveeus
 \scriptfont\eusfam=\eighteus
 \scriptscriptfont\eusfam=\sixeus
\def\fs#1{\fam\eusfam\relax#1}
\def\onsw{{\fs{W}}}
\font\cyr=wncyr10 scaled 1200	
\font\scyr=wncyi10 scaled 1200	

\def\sprod#1#2{{\left\langle{#1},{#2}\right\rangle}}	
\def\txt#1{\quad\mbox{#1}\quad}
\def\d{\mbox{\boldmath$\delta$}}	
\def\m{\mbox{\boldmath$\mu$}}		
\def\eps{\varepsilon}
\def\be{\begin{eqnarray*}}
\def\ee{\end{eqnarray*}}

\newenvironment{proof}{\begin{trivlist}\item[\hskip
	\labelsep{\sc Proof:}]}{
		$\Box$\end{trivlist}}
\newtheorem{corollary}{Corollary}

\newtheorem{theorem}{Theorem}
\newtheorem{lemma}{Lemma}

\begin{document}

\author{J\"org Wenzel}
\date{}
\title{Mean Convergence of Vector--valued\\Walsh Series}
\maketitle

\parskip1.5ex plus .5ex minus .5ex

\begin{abstract}
Given any Banach space $X$, let $L_2^X$ denote the Banach space of all
measu\-rable functions $f:[0,1]\to X$ for which
\[
	\|f\|_2:=\left(\int_0^1\|f(t)\|^2\,dt\right)^{1/2}
\]
is finite. We show that $X$ is a UMD--space (see \cite{BUR:1986}) if and
only if
\[
	\lim_n\|f-S_n(f)\|_2=0\txt{for all $f\in L_2^X$,}
\]
where
\[
	S_n(f):=\sum_{i=0}^{n-1}\sprod{f}{w_i}w_i
\]
is the $n$--th partial sum associated with the Walsh system $(w_i)$.
\end{abstract}

\section{Introduction}
There are several enumerations of the system of Walsh functions. Therefore we
first give the appropriate definition.

For $i\ge1$, the Rademacher functions $(r_i)$ are defined as follows
\[
	r_1(t):=\left\{\begin{array}{cl}
		+1&\txt{for $t\in[0,\frac{1}{2})+\Z$}\\
		-1&\txt{for $t\in[\frac{1}{2},1)+\Z$}
		\end{array}\right.\txt{and}r_{i+1}(t):=r_1(2^it).
\]
Let $n\in\N$. Then $n$ has a unique representation of the form
\[
	n=\sum_{i=0}^\infty n_i2^i,
\]
with $n_i\in\{0,1\}$. Note that in fact only a finite number of the $n_i$
are different from zero. We let
\[
	w_n(t):=\prod_{i=0}^\infty r_{i+1}(t)^{n_i}.
\]
Again the formally infinite product on the righthand side is finite, hence
$w_n$ is well defined.

For $f\in L_2^X$, we denote by $S_n(f)$ the $n$--th partial sum of the Walsh
series of $f$,
\[
	S_n(f):=\sum_{i=0}^{n-1}\sprod{f}{w_i}w_i,
\]
\begin{samepage}
where
\[
	\sprod{f}{w_i}:=\int\limits_0^1f(t)w_i(t)\,dt.
\]
\end{samepage}

Let $X$ and $Y$ be Banach spaces. For $T:X\to Y$, the ideal norm
$\d(T|\onsw_n,\onsw_n)$ is defined to be the least constant $c\ge0$ such that
for all $f\in L_2^X$
\[
	\left\|TS_n(f)\right\|_2\le c\|f\|_2.
\]

Note that
\begin{equation}
\label{norm1}
	\d(T|\onsw_{2^p},\onsw_{2^p})=\|T\|
\end{equation}
for all operators $T$ and $p\in\N$; see e.~g.~\cite{GES}. In order to get a
non--decreasing sequence of ideal norms, we let
\[
	\d^{\max}(T|\onsw_n,\onsw_n):=\max_{1\le k\le n}\d(T|\onsw_k,\onsw_k).
\]
For a more general treatment of ideal norms associated with orthogonal systems
we refer to \cite{PIW}, from where the above notation is adopted.

For $k=1,\dots,2^n$, let
\[
	\Delta_k^{(n)}:=\left[\frac{k-1}{2^n},\frac{k}{2^n}\right)
\]
be the $k$--th dyadic intervall of order $n$.

A dyadic martingale is a martingale $(M_0,M_1,\dots)$ relative to the dyadic
filtration ${\cal F}=({\cal F}_n)$, where ${\cal F}_n$ is generated by $\{
\Delta_k^{(n)}:k=1,\dots,2^n\}$. If $(M_0,M_1,\dots)$ is an $X$--valued dyadic
martingale, then there exist elements $x_j\in X$ such that
\begin{equation}\label{4}
	M_i=\sum_{j=0}^{2^i-1}x_jh_j,
\end{equation}
\begin{samepage}
where $h_j$ denotes the $j$--th Haar function
\[
	h_0\equiv1
\]
\[
\renewcommand{\arraystretch}{1.2}
	h_j(t):=\left\{\begin{array}{cl}
		+2^{(p-1)/2}&\txt{for $t\in\Delta_{2m+1}^{(p)}$,}\\
		-2^{(p-1)/2}&\txt{for $t\in\Delta_{2m+2}^{(p)}$,}\\
		0&\txt{otherwise,}\end{array}\right.
\]
and $j=2^{p-1}+m$, $m=0,\dots,2^{p-1}-1$.
\end{samepage}

As usual, we let
\[
	dM_i:=M_{i+1}-M_i.
\]

Given $p\in\{1,2,\dots\}$, let $\m_p(T)$ denote the least constant $c\ge0$ such
that for all $X$--valued dyadic martingales $(M_0,M_1,\dots,M_p)$ and for all
sequences $\eps_0,\dots,\eps_{p-1}$ of signs $\pm1$ we have
\[
	\left\|\sum_{i=0}^{p-1}\eps_iTdM_i\right\|_2\le c\|M_p\|_2.
\]
We write $\m_p(X)$ instead of $\m_p(I_X)$, where $I_X$ denotes the identity map
of the Banach space $X$.

Note that for all $T:X\to Y$
\begin{equation}
\label{monotonicity}
	\m_{p-1}(T)\le\m_p(T).
\end{equation}
Choosing $M_p:=M_{p-1}$ in the defining inequality of $\m_p(T)$,
we get $dM_{p-1}=0$ and hence
\[
	\left\|\sum_{i=0}^{p-2}\eps_iTdM_i\right\|_2\le\m_p(T)
	\|M_{p-1}\|_2,
\]
which proves the desired inequality.

With the above notation we can prove the following result.
\renewcommand{\thetheorem}{}
\begin{theorem}\label{theorem1}
For all operators $T:X\to Y$ and $p\in\N$, we have
\[
	\d^{\max}(T|\onsw_{2^p},\onsw_{2^p})\le\m_p(T)\le2\d^{\max}(T|\onsw_{
	2^p},\onsw_{2^p}).
\]
\end{theorem}

By definition a Banach space $X$ has the UMD--property if there exists a
constant $c\ge0$ such that
\[
	\left\|\sum_{i=0}^n\eps_idM_i\right\|_2\le c\left\|
	\sum_{i=0}^ndM_i\right\|_2
\]
for all martingales $(M_0,M_1,\dots)$ with values in $X$ and all $n\in\N$. This
is equivalent to the boundedness of the sequence $\m_p(X)$; see
\cite{BUR:1986}.

Thus the theorem gives a characterization of UMD--spaces by the mean
con\-vergence of $X$--valued Walsh series.

\section{Preliminaries}

Let $s,t\in[0,1]$. Then $s$ and $t$ have unique representations $s=\sum_{j=0}
^\infty s_j2^{-j-1}$ and $t=\sum_{j=0}^\infty t_j2^{-j-1}$, respectively,
supposed we choose them to be finite if possible. By $s\oplus t$ we denote
the dyadic sum of $s$ and $t$,
\[
	s\oplus t:=\sum_{j=0}^\infty |t_j-s_j|2^{-j-1}.
\]
Then
\begin{equation}\label{6}
	\int\limits_0^1f(s)\,ds=\int\limits_0^1f(s\oplus t)\,ds
\end{equation}
for all $f\in L_1$ and $t\in [0,1]$. Moreover
\begin{equation}\label{7}
	w_n(s\oplus t)=w_n(s)w_n(t).
\end{equation}

For $n\ge1$, let
\[
	D_n(t):=\sum_{i=0}^{n-1}w_i(t)
\]
be the $n$--th Dirichlet kernel associated with the Walsh functions.

We have
\begin{equation}\label{Sn}
	S_n(f)(t)=\int\limits_0^1f(s)D_n(s\oplus t)\,ds.
\end{equation}
For $n\ge1$, let $0\le k_1<k_2<\dots <k_s$ be defined by
\begin{equation}\label{n}
	n=\sum_{l=1}^s2^{k_l}.
\end{equation}

We will use the following result from \cite[Theorem 8, p.~28]{SWS}.
\begin{lemma}\label{lem5}
\[
	D_n=w_n\sum_{i\in\{k_1,\dots,k_s\}}(D_{2^{i+1}}-D_{2^i}).
\]
\end{lemma}

\section{Proof of the theorem}

For $n$ as in (\ref{n}), we have by (\ref{7}), (\ref{Sn}) and lemma \ref{lem5}
\begin{eqnarray}
\label{TSn}
	\left\|TS_n(f)\right\|_2&=&\left\|\sum_{i\in\{k_1,\dots,k_s\}}\big(TS
	_{2^{i+1}}(fw_n)-TS_{2^i}(fw_n)\big)\right\|_2\nonumber\\&\le&
	\frac{1}{2}\left\|\sum_{i=0}^{k_s}\big(TS_{2^{i+1}}(fw_n)-TS_{2^i}(f
	w_n)\big)\right\|_2+\nonumber\\& &\mbox{}+
	\frac{1}{2}\left\|\sum_{i=0}^{k_s}\eps_i\big(TS_{2^{i+1}}
	(fw_n)-TS_{2^i}(fw_n)\big)\right\|_2,
\end{eqnarray}
where $\eps_i$ is defined by
\[
	\eps_i:=\left\{\begin{array}{cl}
		+1&\txt{if}i\in\{k_1,\dots,k_s\}\\
		-1&\txt{if}i\notin\{k_1,\dots,k_s\}
		\end{array}\right..
\]
Note that $M_i:=S_{2^i}(fw_n)$ form a dyadic martingale of the form (\ref{4}),
since the linear span of the Walsh functions $w_0,\dots,w_{2^p-1}$ coincides
with the linear span of the Haar functions $h_0,\dots,h_{2^p-1}$; see
\cite{PAL:1932}. Hence we have
\be
	\left\|\sum_{i=0}^{k_s}\eps_i
	\big(TS_{2^{i+1}}(fw_n)-TS_{2^i}(Tfw_n)\big)\right\|_2&=&
	\left\|\sum_{i=0}^{k_s}\eps_iTdM_i\right\|_2\\
	&\le&\m_{k_s+1}(T)\left\|S_{2^{k_s+1}}(fw_n)\right\|_2.
\ee
The same argument applied with $\eps_i=+1$ for all $i=0,\dots,k_s$ yields
\be
	\left\|\sum_{i=0}^{k_s}\big(TS_{2^{i+1}}(fw_n)-
	TS_{2^i}(fw_n)\big)\right\|_2&=&\left\|\sum_{i=0}^{k_s}TdM_i
	\right\|_2\\&\le&\m_{k_s+1}(T)\left\|S_{2^{k_s+1}}(fw_n)\right\|_2.
\ee
Therefore we obtain from (\ref{TSn}) that
\[
	\left\|TS_n(f)\right\|_2\le\m_{k_s+1}(T)\left\|S_{2^{k_s+1}}(fw_n)
	\right\|_2.
\]
If $n<2^p$, then $k_s+1\le p$ and it follows from (\ref{norm1}) and
(\ref{monotonicity}) that
\[
	\|TS_n(f)\|_2\le\m_{k_s+1}(T)\|fw_n\|_2\le\m_p(T)\|f\|_2.
\]
If $n=2^p$, then again by (\ref{norm1})
\[
	\|TS_n(f)\|_2\le\|T\|\|S_n(f)\|_2\le\m_p(T)\|f\|_2.
\]
Consequently
\[
	\|TS_n(f)\|_2\le\m_p(T)\|f\|_2
\]
for all $1\le n\le2^p$ and hence
\[
	\d^{\max}(T|\onsw_{2^p},\onsw_{2^p})\le\m_p(T).
\]
This proves the lefthand inequality of the theorem.

To check the righthand inequality we use the following fact.

\begin{lemma}\label{newlem}
Let $I\subseteq\{0,\dots,p-1\}$ and let $n$ be defined by $n:=\sum_{i\in I}
2^i<2^p$. Then we have
\[
	\left\|\sum_{i\in I}TdM_i\right\|_2\le\d(T|\onsw_n,\onsw_n)\|M_p\|_2
\]
for all martingales $(M_0,\dots,M_p)$ of the form {\rm (\ref{4})}.
\end{lemma}

\begin{proof}
We write $M_i$ in the form
\[
	M_i=\sum_{j=0}^{2^i-1}x_jw_j,
\]
where
\[
	x_j:=\int\limits_0^1M_p(t)w_j(t)\,dt\in X.
\]
Then, by lemma \ref{lem5},
\be
	\left\|\sum_{i\in I}TdM_i\right\|_2&=&\left\|\sum_{i\in I}
		\left(\sum_{j=2^i}^{2^{i+1}-1}Tx_jw_j\right)\right\|_2\\
	&=&\left\|\sum_{i\in I}\big(TS_{2^{i+1}}(M_p)-
		TS_{2^i}(M_p)\big)\right\|_2\\
	&=&\left\|\sum_{i\in I}\left(\int\limits_0^1TM_p(s)
		\left(D_{2^{i+1}}(s\oplus t)-D_{2^i}(s\oplus t)\right)\,ds
		\right)\right\|_2\\
	&=&\left\|\int\limits_0^1TM_p(s)w_n(s\oplus t)D_n(s\oplus t)\,ds
	\right\|_2\\
	&=&\|TS_n(M_pw_n)\|_2\\
	&\le&\d(T|\onsw_n,\onsw_n)\|M_pw_n\|_2\\
	&=&\d(T|\onsw_n,\onsw_n)\|M_p\|_2.
\ee
\end{proof}

We are now able to complete the proof of the theorem.
Let a sequence $\eps_0,\dots,\eps_{p-1}$ of signs $\pm1$ be
given. Define $n$ and $m$ by
\[
	n:=\sum_{\{i:\eps_i=+1\}}2^i\txt{and}
	m:=\sum_{\{i:\eps_i=-1\}}2^i.
\]
Then we get from lemma \ref{newlem} that
\be
	\left\|\sum_{i=0}^{p-1}\eps_iTdM_i\right\|_2&\le&
	\left\|\sum_{\{i:\eps_i=+1\}}TdM_i\right\|_2+
	\left\|\sum_{\{i:\eps_i=-1\}}TdM_i\right\|_2\\
	&\le&\d(T|\onsw_n,\onsw_n)\|M_p\|_2+\d(T|\onsw_m,\onsw_m)\|M_p\|_2\\
	&\le&2\d^{\max}(T|\onsw_{2^p},\onsw_{2^p})\|M_p\|_2.
\ee
Since this holds for all sequences $(\eps_i)$, we have
\[
	\m_p(T)\le 2\d^{\max}(T|\onsw_{2^p},\onsw_{2^p}),
\]
which is the desired righthand inequality.

\section{Some consequences}

\begin{corollary}\label{corollary2}
The following conditions are equivalent.
\begin{enumerate}
\item
$\|f-S_n(f)\|_2\to 0$ for all $f\in L_2^X$.
\item
$X$ has the UMD--property.
\end{enumerate}
\end{corollary}

\begin{proof}
If $X$ has the UMD--property, then by the theorem
\[
	\d^{\max}(X|\onsw_n,\onsw_n)\le c
\]
for all $n\in\N$. Since the Walsh functions form a complete orthonormal system
in $L_2[0,1]$, we can find a linear combination $\sum_{k=0}^Nx_kw_k\in
L_2^X$ with
\[
	\left\|f-\sum_{k=0}^Nx_kw_k\right\|_2\le\eps.
\]
Then, for $n\ge N$,
\begin{eqnarray*}
	\left\|f-S_n(f)\right\|_2&\le &\left\|f-\sum_{k=0}^Nx_kw_k\right\|_2+
	\left\|S_n(\sum_{k=0}^Nx_kw_k-f)\right\|_2\\
	&\le&\eps+\d^{\max}(X|\onsw_n,\onsw_n)\eps
	\le (1+c)\eps.
\end{eqnarray*}
If on the other hand
\[
	\|f-S_n(f)\|_2\to0\txt{for all $f\in L_2^X$,}
\]
then, by the uniform boundedness theorem, we get
\[
	\|S_n(f)\|_2\le c\|f\|_2
\]
and hence $X$ is a UMD--space.
\end{proof}

As a further easy application of the theorem we get the order of growth of
$\m_p(X)$ for $X=L_1[0,1]$.

To this end, let
\begin{equation}\label{lebesgue}
	L_n:=\int\limits_0^1|D_n(t)|\,dt
\end{equation}
be the Lebesgue constants associated with the Walsh system. We also consider
\[
	L_n^{\max}:=\max_{k\le n}L_k.
\]

\begin{corollary}\label{theorem3}
\[
	\frac{p}{8}\le\m_p(L_1[0,1])\le2p.
\]
\end{corollary}

\begin{proof}
For $X=L_1[0,1]$ and $f\in L_1[0,1]$, we define $F\in L_2^X$ by
\[
	F(t):=f(t\oplus\cdot).
\]
Then, by (\ref{6}) and (\ref{7}), we have
\be
	\sprod{F}{w_j}&=&\int\limits_0^1F(t)w_j(t)\,dt=\int\limits_0^1f(t\oplus
	\cdot)w_j(t)\,dt\\&=&\int\limits_0^1f(t)w_j(t)w_j(\cdot)\,dt=\sprod{f}{
	w_j}w_j.
\ee
Hence
\[
	\sum_{j=0}^{n-1}\sprod{F}{w_j}w_j(t)=\sum_{j=0}^{n-1}\sprod{f}{w_j}
	w_j(t\oplus\cdot)
\]
Furthermore 
\be
	\|F\|_2^2& = &\int\limits_0^1\|F(t)\|_1^2dt=\int\limits_0^1\left(\int
	\limits_0^1|f(t\oplus s)|\,ds\right)^2dt\\&=&\int\limits_0^1\left(\int
	\limits_0^1|f(s)|\,ds\right)^2dt=\|f\|_1^2.
\ee
Similarly
\[
	\left\|\sum_{j=0}^{n-1}\sprod{F}{w_j}w_j\right\|_2=
	\left\|\sum_{j=0}^{n-1}\sprod{f}{w_j}w_j\right\|_1.
\]
If we now choose $f$ to be the characteristic function of the intervall
$[0,2^{-p}]$, then for $k\le n<2^p$
\[
	\sprod{f}{w_k}=2^{-p}\txt{and}\|f\|_1=2^{-p}.
\]
Consequently
\be
	\d(L_1|\onsw_n,\onsw_n)&\ge&\frac{\left\|\sum_{j=0}^{n-1}\sprod{F}{w_j}
	w_j\right\|_2}{\|F\|_2}\\&=&L_n,
\ee
where $L_n$ denotes the Lebesgue constant as defined in (\ref{lebesgue}).

Since
\[
	\d(X|\onsw_n,\onsw_n)\le L_n
\]
for all Banach spaces $X$, we have
\[
	\d(L_1|\onsw_n,\onsw_n)=L_n.
\]

This proves corollary \ref{theorem3} by taking into account our theorem
and the following result from \cite[Theorem 9, p.~34]{SWS}.
\begin{lemma}\label{lem6}
\[
	\frac{p}{8}\le L_{2^p}^{\max}\le p.
\]
\end{lemma}
\end{proof}

\noindent
Mathematische Fakult\"at\\
Friedrich--Schiller--Universit\"at Jena\\
Universit\"atshochhaus 17. OG\\
O--6902 Jena\\
Germany\\
e--mail: wenzel@mathematik.uni--jena.de.dbp (X400)

\end{document}